\documentclass[review,3p,10pt]{elsarticle}

\biboptions{sort&compress}

\usepackage{amssymb}
\usepackage{amsmath}
\usepackage{cases}
\usepackage{epstopdf}
\usepackage{fixltx2e}
\usepackage{mathtools}
\usepackage{array}
\usepackage{bm}
\usepackage{multirow}
\usepackage{color}
\usepackage{algorithm}
\usepackage{algorithmic}
\usepackage{amsthm}
\usepackage{mathrsfs}
\usepackage{amssymb}

\newtheorem{assumption}{Assumption}

\newtheorem{remark}{Remark}
\usepackage{lineno}

\usepackage{url}

\newcommand{\tabincell}[2]{\renewcommand\arraystretch{0.8}\begin{tabular}{@{}#1@{}}#2\end{tabular}}

\journal{arXiv.org}

\begin{document}
\begin{frontmatter}



\title{Mitigating Renewable-Induced Risks for Green and Conventional Ammonia Producers through Coordinated Production and Futures Trading}

\author[label1]{Huayan~Geng}
\author[label2]{Yangjun~Zeng\corref{cor1}}
\ead{zengyangjun@stu.scu.edu.cn}
\author[label2]{Yiwei~Qiu\corref{cor2}}
\ead{ywqiu@scu.edu.cn}

\address[label1]{School of Statistics and Data Science,	Southwestern University of Finance and Economics, Chengdu 611130, China }
\address[label2]{College of Electrical Engineering, Sichuan University, Chengdu, 610065, China}
\cortext[cor1]{Corresponding author}
\cortext[cor2]{Corresponding author}

\begin{abstract}
Renewable power-to-ammonia (ReP2A), which uses hydrogen produced from renewable electricity as feedstock, is a promising pathway for decarbonizing the energy, transportation, and chemical sectors. However, variability in renewable generation causes fluctuations in hydrogen supply and ammonia production, leading to revenue instability for both ReP2A producers and  conventional fossil-based gray ammonia (GA) producers in the market. Existing studies mainly rely on engineering measures, such as production scheduling, to manage this risk, but their effectiveness is constrained by physical system limits. To address this challenge, this paper proposes a financial instrument termed \emph{renewable ammonia futures}  and integrates it with production decisions to hedge ammonia output risk. Production and trading models are developed for both ReP2A and GA producers, with conditional value-at-risk (CVaR) used to represent risk preferences under uncertainty. A game-theoretic framework is established in which the two producers interact in coupled ammonia spot and futures markets, and a Nash bargaining mechanism coordinates their production and trading strategies. Case studies based on a real-world system show that introducing renewable ammonia futures increases the CVaR utilities of ReP2A and GA producers by 5.103\% and 10.14\%, respectively, improving profit stability under renewable uncertainty. Sensitivity analysis further confirms the effectiveness of the mechanism under different levels of renewable variability and capacity configurations.
\end{abstract}


\end{frontmatter}

\section{Introduction}
\label{sec:intro}

\subsection{Motivation}
\label{sec:motivation}

Amid the global energy transition, renewable power-to-hydrogen (ReP2H) has been widely recognized as a key pathway for integrating renewable electricity \cite{lima2026dynamic} and decarbonizing hard-to-abate sectors such as transport and chemicals \cite{ajanovic2024on,jovan2022utilization,li2022coplanning}. As the largest downstream product of hydrogen, ammonia is increasingly considered an energy carrier and industrial feedstock for low-carbon systems \cite{li2025development,li2025redesigning}. Consequently, renewable power-to-ammonia (ReP2A) projects are expanding worldwide. In China alone, more than 90 green hydrogen and renewable ammonia (RA) projects have been planned, with a combined capacity of approximately 15 million tons per year \cite{lin2025review,daan2025green,etog2025shenneng,jilin2025bridging}. Large-scale projects are also emerging in Europe and other regions \cite{narciso2025design}.

In ReP2A systems, renewable electricity drives water electrolysis to produce hydrogen, which is subsequently converted to ammonia through an improved Haber-Bosch process \cite{li2025redesigning}. However, the variability of renewables leads to fluctuations in hydrogen supply and ammonia output. Historical data show that annual wind and solar generation can vary by 30\%-40\% between high- and low-generation years, with even larger variations at monthly timescales \cite{yu2023optimal,yu2024optimal}.

Such variability introduces economic risks in ammonia markets. For ReP2A producers, uncertain production results in unstable sales and revenue. For conventional fossil-based gray ammonia (GA) producers, renewable variability does not directly affect production but increases market price volatility through fluctuations in RA supply \cite{wu2023multitimescale,cai2023green}. Maintaining stable profitability under uncertain RA output therefore becomes a key challenge for both producers.

Existing research mainly addresses this problem through engineering solutions. 
Typical approaches include:
1) installing electricity, hydrogen, or ammonia storage to smooth supply and output over time \cite{wen2022technoeconomic,du2025comprehensive,rosbo2023flexible}; and
2) optimizing scheduling and operational control under renewable variability \cite{wu2024research,wang2025optimal,shi2025optimal}.
However, these measures are constrained by physical and engineering limits. Electrochemical energy and hydrogen storage typically provide balancing only at short timescales because of cost and safety considerations \cite{wen2022technoeconomic,du2025comprehensive}. Although ammonia can be stored in liquid form, storage tanks are classified as hazardous facilities and are subject to strict safety regulations, which limits their capacity. Therefore, shifting production across monthly or seasonal timescales remains difficult.

To complement engineering measures, this study introduces a financial instrument termed \emph{renewable ammonia futures}, inspired by hedging mechanisms in renewable electricity markets. Electricity futures have long been used to manage price and production risks \cite{bessembinder2002equilibrium,botterud2010relationship}. In countries such as Germany, renewable futures contracts allow wind power producers and conventional generators to agree on delivery quantities and prices for a future period. Such contracts help renewable producers hedge uncertainty while enabling conventional generators to manage demand fluctuations caused by volatile renewable output \cite{gersema2017equilibrium,benth2009dynamic,liu2019pricing}.
These practices motivate the exploration of similar instruments for ammonia markets. However, ammonia production involves different physical processes and market structures from electricity systems, and the effectiveness of futures-based risk management remains unclear.

Building on these observations, this study proposes a renewable ammonia futures mechanism inspired by renewable electricity futures \cite{gersema2017equilibrium,benth2009dynamic,liu2019pricing} and our previous work \cite{cai2023green}. Considering both the physical production characteristics of ReP2A and its interactions with the GA producer, a coupled spot-futures market framework is developed. Through equilibrium analysis and bargaining mechanisms, this study investigates how engineering measures and financial instruments can jointly mitigate production risks.


\subsection{Literature review}
\label{sec:review}

Extensive research has examined how ReP2A systems can accommodate renewable variability. Existing studies mainly focus on three areas: improving process flexibility, optimizing storage capacities, and developing operational scheduling strategies.

\subsubsection{Enhancing flexibility of hydrogen and ammonia production processes}

Many studies aim to improve the operational flexibility of hydrogen and ammonia production while maintaining safe operation. For hydrogen production, Xia et al. \cite{xia2024efficiency} improved the efficiency and load range of alkaline electrolyzers (AELs) by developing pulse rectifiers and coordinated temperature-pressure controller. Sha et al. \cite{sha2026cascaded} proposed a cascaded process integrating multiple stacks with lye-gas separators to improve current efficiency. Hu et al. \cite{hu2025analysis} analyzed electrochemical characteristics and impurity accumulation under varying operating conditions and identified safe operating boundaries for AELs.
For ammonia synthesis, Fahr et al. \cite{fahr2025dynamic} developed dynamic models of the Haber-Bosch process and optimized reactor configurations to expand the load range. Subsequent work \cite{fahr2023design,verleysen2023how} further improved reactor thermodynamic design to enable low-load operation. Ji et al. \cite{ji2024multistable} and Zhang et al. \cite{zhang2024how} proposed multi-steady-state ammonia synthesis processes suitable for flexible ReP2A operation.

While these studies improves operational flexibility, they mainly enable systems to follow renewable fluctuations rather than reduce the economic risks associated with production variability.

\subsubsection{Sizing of electrical energy, hydrogen, and ammonia storage}

Energy and product storage are widely used to buffer renewable variability. Battery energy storage systems (BESSs) provide fast regulation from seconds to hours, while hydrogen storage supports longer-term balancing. Many studies therefore optimize storage configurations to address variability across multiple timescales \cite{rosbo2023flexible,yu2023optimal,yu2024optimal,zhou2025flexible,lan2024intelligent,zeng2025planning}.

In principle, sufficiently large storage capacities could stabilize hydrogen supply and ammonia production. In practice, however, storage capacity is limited by cost and safety considerations \cite{wen2022technoeconomic,du2025comprehensive}. For example, installing a BESS equivalent to ten minutes of electrolyzer rated power can increase the levelized cost of ammonia (LCOA) by about 200 CNY/t \cite{zhu2026exploring}. Hydrogen storage also introduces additional infrastructure costs and safety requirements.
When hydrogen storage capacity reaches 100,000 Nm$^3$, safety risks increase significantly \cite{le2023safety}. According to the Chinese national standard \emph{GB/T 29729-2022 Safety Technical Specification for Pressure Gaseous Hydrogen Storage Systems}, such facilities must satisfy strict requirements regarding siting and safety distance. Even for a plant producing 200,000 t of ammonia annually, a hydrogen storage capacity of 500,000 Nm$^3$ can sustain full-load synthesis for only several hours.

Liquid ammonia storage provides longer-term buffering but is also constrained by safety regulations because ammonia tanks are classified as hazardous facilities. Thus, storage alone cannot mitigate risks arising from inter-month renewable variability \cite{yu2023optimal,yu2024optimal,wu2023multitimescale,wu2024dispatchable}.

\subsubsection{Optimal scheduling and control}

Another research direction focuses on coordinated scheduling of renewable generation, hydrogen production, ammonia synthesis, and storage. Qiu et al. \cite{qiu2023extended} and Zeng et al. \cite{zeng2025scheduling,zeng2025optimal,zeng2026harmonic} coordinated electrolyzers and rectifiers while considering electrochemical safety and grid constraints. Rosbo et al. \cite{rosbo2023flexible} developed dynamic models of ammonia production systems and integrated stability analysis with scheduling. Kong et al. \cite{kong2024nonlinear} proposed a model predictive controller for the Haber-Bosch process. At the system level, Shi et al. \cite{shi2025optimal} developed a chance-constrained scheduling model for ReP2A systems. Wu et al. \cite{wu2024research} proposed a capacity planning and scheduling framework with annual profit maximization. Wang \cite{wang2025optimal} formulated a stochastic two-stage mixed-integer programming model for power-to-ammonia scheduling. Wu et al. \cite{wu2023multitimescale,wu2024dispatchable} further examined multi-timescale scheduling strategies spanning annual to intra-day horizons.

Although these approaches improve operational adaptability, they remain constrained by physical limits and cannot fully eliminate revenue risks caused by renewable variability.

\vspace{6pt}

In electricity markets, financial instruments such as renewable futures contracts have been introduced to hedge risks associated with variable generation \cite{gersema2017equilibrium,benth2009dynamic,liu2019pricing}.
For example, Benth et al. \cite{benth2009dynamic} analyzed the impact of seasonal wind output on wind futures pricing. Gersema et al. \cite{gersema2017equilibrium} and Liu et al. \cite{liu2019pricing} further examined pricing and hedging effects using equilibrium models.
Inspired by these mechanisms, our previous work \cite{cai2023green} explored their potential application to ammonia markets. However, \cite{cai2023green} did not incorporate engineering constraints or the production decisions of ReP2A and GA producers.

To address these limitations, this study develops an integrated framework that combines engineering measures with renewable ammonia futures mechanisms to mitigate production risks. The contributions are briefed in Section \ref{sec:contribution}.

\subsection{Contributions of this work}
\label{sec:contribution}

This study develops an integrated framework combining engineering measures and financial instruments to mitigate risks caused by renewable variability in ammonia production. The main contributions are:

\begin{enumerate}
	\item A financial instrument termed \emph{renewable ammonia futures} is proposed, inspired by renewable electricity futures markets. The mechanism complements engineering measures by enabling producers to hedge risks associated with RA output fluctuations.
	
	\item A game-theoretic framework incorporating conditional value-at-risk (CVaR) is developed to model the interactions between ReP2A and fossil-based GA producers in coupled spot and futures markets. A Nash bargaining strategy coordinates production and trading decisions while improving the utilities of both participants.
	
	\item Case studies based on a real-world system show that integrating futures trading increases the utilities of ReP2A and GA producers by 5.103\% and 10.14\%, respectively. Sensitivity analyses further examine how renewable uncertainty and capacity configurations affect the performance of the mechanism.
\end{enumerate}

The remainder of this paper is organized as follows. Section \ref{sec:problem} introduces the system and market structures. Section \ref{sec:futures} presents the renewable ammonia futures mechanism. Section \ref{sec:optimization} develops the production and trading decision models. Section \ref{sec:game} describes the spot-futures interaction model and solution method. Section \ref{sec:case} presents the case studies, and Section \ref{sec:conclusion} concludes the paper.

\section{Problem description}
\label{sec:problem}
	
\subsection{Overview of ReP2A and GA systems}
\label{sec:sys}

\begin{figure}[t]
	\centering
	\includegraphics[scale=1.1]{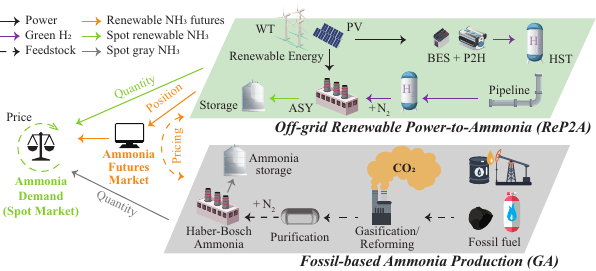}\vspace{-6pt}
	\caption{Illustration of ReP2A and GA systems and their interactions with ammonia futures and spot markets.}
	\label{fig:sys}
\end{figure}

The structures of the ReP2A and GA systems are illustrated in Fig. \ref{fig:sys}. A typical ReP2A system includes renewable generation, BESS, water electrolysis for hydrogen production, air separation for nitrogen production, hydrogen compression, ammonia synthesis, and liquid ammonia storage \cite{narciso2025design}. In practice, these units are usually invested in and operated by a single entity as an integrated facility.

In the ReP2A process, renewable electricity drives electrolytic hydrogen production. The hydrogen is then combined with nitrogen from air separation and fed to the ammonia synthesis unit, where ammonia is produced through an improved Haber-Bosch process \cite{ji2024multistable,zhang2024how}. BESS, hydrogen storage, and ammonia storage provide buffering at different time scales. Due to policy and regulatory constraints, ReP2A systems are typically operated off-grid and cannot trade electricity with the power grid \cite{zeng2025planning,yu2023optimal,Campion2024potential}.

The GA system follows the conventional Haber-Bosch process. It includes syngas production from coal or natural gas, air separation, ammonia synthesis, and liquid ammonia storage. Because fossil fuels can be stored reliably, GA production is largely independent of renewable power variability and can operate steadily. Output levels may also be adjusted in response to market conditions.

\subsection{Ammonia market structure and assumptions}
\label{sec:assumption}

RA is currently traded in two main markets, i.e., the traditional chemical market (e.g., feedstock for urea and ammonium salts) and the emerging maritime fuel market \cite{li2025redesigning}. In the shipping sector, RA can avoid carbon taxes and obtain a green premium, allowing it to be sold at prices higher than fossil-based GA. In contrast, the chemical market does not yet differentiate between RA and GA. They are therefore treated as homogeneous products and traded at the same price. This study focuses on the chemical market. Without loss of generality, the following assumptions are made.

\begin{assumption}
	\label{assum:1}
	RA and GA are homogeneous products and are sold at the same price.
\end{assumption}

\begin{assumption}
	\label{assum:2}
	Because ammonia is a hazardous chemical with limited transportation range, a local market with one ReP2A producer and one GA producer is considered.
\end{assumption}

\begin{assumption}
	\label{assum:3}
	Ammonia trading occurs monthly. 
	The study horizon is one year with 12 trading periods.
\end{assumption}

\begin{assumption}
	\label{assum:4}
	Due to strict regulatory limits on hazardous chemical storage \cite{wu2023multitimescale}, ammonia produced in each month is assumed to be sold within the same month.
\end{assumption}

\begin{assumption}
	\label{assum:5}
	Ammonia demand follows a commonly used price elasticity model \cite{zeng2025planning}:
	\begin{align}
	\rho _t^{{\rm{am}}} = {\rho ^{{\rm{max}}}} - (M_t^{{\rm{ga,sell}}} + M_t^{{\rm{ra,sell}}})/{k^{{\rm{am}}}}, \label{eq:elasticity}
	\end{align}
	\noindent
	where $t=1, \ldots, 12$ denotes the trading period; $\rho_t^{\rm{am}}$ is the ammonia price; $M_t^{{\rm{ga,sell}}}$ and $M_t^{{\rm{ra,sell}}}$ are the quantities supplied by the GA and ReP2A producers; ${\rho ^{{\rm{max}}}}$ is the intercept price when supply is zero; and $k^{\rm{am}}$ is the price elasticity coefficient.
\end{assumption}

\subsection{Risk metrics for production and trading decisions}
\label{sec:risk}

Variability in renewable generation introduces uncertainty in ReP2A production and creates economic risks for both producers. For ReP2A, output fluctuations lead to unstable revenue. For the GA producer, renewable variability does not directly affect production but alters RA supply and thus increases ammonia price volatility. Production decisions must therefore account for both expected profit and downside risk.

To quantify risk, this study adopts the conditional value-at-risk (CVaR), which measures expected losses in the tail of the profit distribution by presenting the expected loss exceeding the value-at-risk (VaR) at a given confidence level $\alpha$, and can be incorporated into the objective function to obtain risk-aware production strategies.

Let $L$ denote the loss function under uncertainty (defined as negative profit; its construction is given in Sections \ref{sec:modelrep2a}--\ref{sec:modelga}). For confidence level $\alpha \in [0,1)$, VaR and CVaR are defined as
\begin{gather}
	{\rm{Va}}{{\rm{R}}_\alpha }\left( L \right) = \inf \left\{ {l \in \,{\rm{:}}\,{\rm{Prob}}\left( {L \le l} \right) \ge \alpha } \right\}, \label{eq:var} \\
	{\rm{CVa}}{{\rm{R}}_\alpha }\left( L \right) = \mathbb{E} \left[ {L\left| {L \ge {\rm{Va}}{{\rm{R}}_\alpha }\left( L \right)} \right.} \right], \label{eq:cvar}
\end{gather}

Following \cite{rockafellar2000optimization,rockafellar2002conditional}, (\ref{eq:cvar}) can be reformulated to facilitate convex optimization:
\begin{gather}
	{\rm{CVa}}{{\rm{R}}_\alpha }\left( L \right) = \mathop {\min }\limits_{\theta  \in \mathbb{R}} \left\{ {\theta  + {1 \over {1 - \alpha }}\left[ {{{\left( {L - \theta } \right)}^ + }} \right]} \right\}, \label{eq:cvartrans}
\end{gather}
where $(\cdot)^+ \triangleq \max(\cdot,0)$; $\theta$ approximates ${\rm{Va}}{{\rm{R}}_\alpha}\left( L \right)$.

When renewable uncertainty is represented by a scenario set $\Omega={1,\ldots,N^\text{s}}$, (\ref{eq:cvartrans}) can be linearized as
\begin{gather}
	{\rm{CVa}}{{\rm{R}}_\alpha }\left( L \right) = \mathop {\min }\limits_{\theta ,{\xi ^\omega }} \left\{ {\theta  + {1 \over {1 - \alpha }}\sum\nolimits_{\omega  \in \Omega } {{p^\omega }{\xi ^\omega }} } \right\}, \label{eq:cvarlin}\\
	\text{s.t.}\ {\xi ^\omega } \ge {L^\omega } - \theta ,\;{\xi ^\omega } \ge 0, \label{eq:cvarcons}
\end{gather}
\noindent
where $p^\omega$, $L^\omega$, and $\xi^\omega$ denote the probability, loss value, and excess loss under scenario $\omega$.

\section{Renewable ammonia futures trading mechanism}
\label{sec:futures}

When physical buffering options such as energy and ammonia storage are limited, intertemporal trading can help redistribute production risk at the market level. Inspired by the German wind power futures mechanism \cite{gersema2017equilibrium}, this study introduces an ammonia futures contract between ReP2A and the GA producer. Under this arrangement, ReP2A transfers the ownership of part of its RA output to the GA producer in advance for a fixed payment, allowing both parties to hedge risks associated with RA production variability. The mechanism operates as follows.

\emph{Step 1:} At the beginning of each trading period $t$, ReP2A negotiates with the GA producer to sell ammonia futures with position $Q_t^{\rm{f}}$ at price $\rho _t^{\rm{f}}$.

\emph{Step 2:} At the end of period $t$, the contract is settled. ReP2A transfers ownership of part of RA production to the GA producer, who then sells the delivered ammonia in the spot market. The delivery quantity can be defined in two ways:

\textbf{a) Mode 1 (wind-futures-type mechanism).}
The position $Q_t^{\rm{f}}$ represents the share of RA production transferred to the GA producer. It is a dimensionless variable between 0 and 1. The delivery quantity is therefore proportional to realized production:
\begin{gather}
  M_t^{\rm{f}} = Q_t^{\rm{f}} \cdot M_t^{{\rm{ra,pro}}}, \label{eq:future1} \\
  0 \le Q_t^{\rm{f}} \le 1, \label{eq:future2}
\end{gather}
\noindent
where $M_t^{\rm{f}}$ denotes the ammonia quantity transferred for contract settlement, measured in tons; and $M_t^{{\rm{ra,pro}}}$ is the RA production in period $t$.

This mechanism provides a natural hedging effect. By selling a share of future production, ReP2A secures part of its revenue and reduces exposure to renewable uncertainty. By acquiring a share of RA output, the GA producer can hedge risks arising from fluctuations in RA supply.

For example, when renewable generation is high, RA output increases and market prices may fall due to demand elasticity. Because the GA producer has already acquired part of the RA output through the futures contract, it can sell more ammonia and partly offset the price decline. Conversely, when renewable generation is low, the delivered RA decreases, but reduced market supply raises prices and compensates for lower sales volume.

\textbf{b) Mode 2 (conventional commodity-futures mechanism).}
In this mode, the futures position is defined as a fixed delivery quantity, i.e., $Q_t^{\rm{f}} = M_t^{\rm{f}}$ measured in tons. The contract specifies a predetermined amount of RA to be transferred at settlement, independent of renewable generation. The delivered ammonia, together with GA production and any remaining RA, is sold in the spot market.

\begin{remark}
	\label{remark:1}
	Under the proposed framework, ammonia futures are traded only between ReP2A and GA producers. Consumers participate only in the spot market.
\end{remark}

\begin{remark}
	\label{remark:2}
	Mode 2 is introduced for comparison. Because the delivery quantity is independent of renewable production, it does not create a negative correlation between profit and production risk. Consequently, no combination of futures price and trading volume can simultaneously improve the utilities of both producers, and the bargaining process fails, as shown in Section \ref{sec:casebargining}.
\end{remark}

For convenience, futures decisions over the entire horizon $t=1,\ldots,T$ are represented as vectors ${\bm{Q}}^{\rm{f}} \triangleq {\left[ {Q_t^{\rm{f}}} \right]_{t = 1, \ldots ,T}}$ and $\bm{\rho}^{\rm{f}} \triangleq {\left[ {\rho_t^{\rm{f}}} \right]_{t = 1, \ldots ,T}}$.

\section{Decision models of ReP2A and GA in the ammonia market under risk}
\label{sec:optimization}

Both ReP2A and GA producers are assumed to be rational and seek to maximize their utilities under risk through production scheduling and futures-spot trading. Their decisions are formulated as mathematical programs using the CVaR metric introduced in Section \ref{sec:risk}.

\subsection{Operational decision of ReP2A}
\label{sec:modelrep2a}

The profit of the ReP2A producer, denoted by ${C^{{\rm{ra}}}}$, is defined as
\begin{gather}
	{C^{{\rm{ra}}}} = \sum\nolimits_{t = 1}^T {\left[ {\left( {M_t^{{\rm{ra,sell}}} - M_t^{\rm{f}}} \right)\rho _t^{{\rm{am}}} + Q_t^{\rm{f}}\rho _t^{\rm{f}} - {c^{{\rm{ra}}}}} \right]}, \label{eq:objrep2a}
\end{gather}
\noindent
where the summation terms represent spot market revenue, futures revenue, and production cost; $M_t^{{\rm{ra,sell}}}$ denotes total RA supplied to the spot market; the quantity $M_t^{{\rm{ra,sell}}}-M_t^{\rm{f}}$ is sold directly by ReP2A, while $M_t^{\rm{f}}$ is transferred to the GA producer through futures settlement; $Q_t^{\rm{f}}\rho_t^{\rm{f}}$ represents futures income; ${c^{{\rm{ra}}}}$ is the production cost. Because renewable electricity is generated within the integrated system, variable costs such as electricity procurement costs are not considered.

ReP2A production must satisfy system-level operational constraints. As discussed in Section \ref{sec:review}, coordinated scheduling of energy storage, electrolysis, and ammonia synthesis allows a stable share of renewable generation to be converted into ammonia \cite{yu2023optimal,yu2024optimal,wu2023multitimescale,wu2024dispatchable}. At the monthly time scale adopted here, device-level dynamics are not modeled explicitly. Instead, the energy conversion and production limits are described as
\begin{gather}
   M_t^{\rm{ra,pro}} \le \eta^{\rm{p2a}} E_t^{\rm{res}}, \label{eq:cons1rep2a} \\
   \underline{M}^{\rm{ra,pro}} \le M_t^{{\rm{ra,pro}}} \le {\overline M^{{\rm{ra,pro}}}}, \label{eq:cons2rep2a} \\
   M_t^{{\rm{ra,sell}}} \le M_t^{{\rm{ra,pro}}}, \label{eq:cons3rep2a}
\end{gather}
\noindent
where $E_t^{\rm{res}}$ denotes renewable electricity available for ammonia production in period $t$; $M_t^{{\rm{ra,pro}}}$ is RA production; $\eta^{\rm{p2a}}$ is the power-to-ammonia conversion factor; and $\underline{M}^{\rm{ra,pro}}$ and $\overline{M}^{\rm{ra,pro}}$ are technical production limits. More detailed ReP2A scheduling model can be found in our previous studies \cite{wu2023multitimescale,zeng2025planning,wu2024dispatchable}.

For compact notation, the operational decision variables of ReP2A are defined as
\begin{gather}
  {\bm{u}^{{\rm{ra}}}} \buildrel \Delta \over = {\left[ {M_t^{{\rm{ra,pro}}},M_t^{{\rm{ra,sell}}}} \right]_{t = 1, \ldots ,T}}.
\end{gather}
\noindent
Because these variables depend on renewable generation realizations, scenario-dependent decisions are denoted by ${[{\bm{u}}_\omega^{{\rm{ra}}}]_{\omega \in \Omega }}$. Futures variables are determined through negotiation and are therefore written separately:
\begin{gather}
	{\bm{v}} \buildrel \Delta \over = {\left[ {Q_t^{\rm{f}},\rho _t^{\rm{f}}} \right]_{t = 1, \ldots ,T}}.
\end{gather}

To account for risk, the objective is to minimize the CVaR of $-C^{{\rm{ra}}}$ while satisfying constraints (\ref{eq:cons1rep2a})--(\ref{eq:cons3rep2a}) for all scenarios. The ReP2A operational decision problem can therefore be formulated:
\begin{gather}
	\mathop {\min }\limits_{{{[{\bm{u}}_\omega ^{{\rm{ra}}}]}_{\omega  \in \Omega }}}
	\ {f^{{\rm{ra}}}}\left( {{{[{\bm{u}}_\omega ^{{\rm{ra}}}]}_{\omega  \in \Omega }},{\bm{v}}} \right) \buildrel \Delta \over = {\rm{CVa}}{{\rm{R}}_\alpha }\left( { - {C^{{\rm{ra}}}}} \right), \label{eq:objcompactrep2a} \\
	\text{s.t.} \ {\bm{h}}_{}^{{\rm{ra}}}\left( {{\bm{u}}_\omega ^{{\rm{ra}}},{\bm{v}},{\bm{\phi}_\omega }} \right) = \bm{0},\forall\omega \in \Omega \;:\;{\bm{\lambda} ^{{\rm{ra}}}} , \label{eq:hconsobjcompactrep2a} \\
	 {\bm{g}}_{}^{{\rm{ra}}}\left( {{\bm{u}}_\omega ^{{\rm{ra}}},{\bm{v}},{\bm{\phi}_\omega  }} \right) \le \bm{0},\forall\omega \in \Omega \;:\;{\bm{\mu} ^{{\rm{ra}}}}, \label{eq:gconsobjcompactrep2a}
\end{gather}
\noindent
where ${f^{{\rm{ra}}}}(\cdot)$ is  referred to as the \emph{utility function};   ${\bm{h}}^{{\rm{ra}}}(\cdot)$ and ${\bm{g}}^{{\rm{ra}}}(\cdot)$ represent all equality and inequality constraints; and $\bm{\phi}_\omega = [E^{\omega}_1,\ldots,E^{\omega}_T]$ denotes the renewable generation vector in scenario $\omega$; ${\bm{\lambda} ^{{\rm{ra}}}}$ and ${\bm{\mu} ^{{\rm{ra}}}}$ are dual variables. After applying the CVaR linearization in (\ref{eq:cvarlin})--(\ref{eq:cvarcons}), the resulting problem becomes a linear program (LP).

\subsection{Operational decision of GA}
\label{sec:modelga}

Similar to ReP2A, the profit function of the GA producer, denoted by ${C^{{\rm{ga}}}}$, is defined as
\begin{gather}
	{C^{{\rm{ga}}}} = \sum\nolimits_{t = 1}^T {\left[ {\left( {M_t^{{\rm{ga,sell}}} + M_t^{\rm{f}}} \right)\rho _t^{{\rm{am}}} - Q_t^{\rm{f}}\rho _t^{\rm{f}} - \left( {c_{\rm{0}}^{{\rm{ga}}} + c_1^{{\rm{ga}}}M_t^{{\rm{ga,pro}}}} \right)} \right]}, \label{eq:objga}
\end{gather}
\noindent
where the summation terms represent spot revenue, payment for futures contracts, and production cost; $M_t^{{\rm{ga,sell}}}$ denotes ammonia sales from GA production, while $M_t^{\rm{f}}$ is the RA obtained through futures settlement. Thus the total spot-market supply becomes $M_t^{{\rm{ga,sell}}}+M_t^{\rm{f}}$. The parameters $c_{\rm{0}}^{{\rm{ga}}}$ and $c_1^{{\rm{ga}}}$ represent fixed and variable production costs.

GA production and sales are subject to
\begin{gather}
 {\underline{M}^{{\rm{ga,pro}}}} \le M_t^{{\rm{ga,pro}}} \le {\overline{M}^{{\rm{ga,pro}}}}, \label{eq:cons1ga} \\
  M_t^{{\rm{ga,sell}}} \le M_t^{{\rm{ga,pro}}}, \label{eq:cons2ga}
\end{gather}
\noindent
where $M_t^{{\rm{ra,pro}}}$, $\underline{M}^{{\rm{ga,pro}}}$, and $\overline{M}^{{\rm{ga,pro}}}$ denote the GA production and technical limits.

The decision variables are written compactly as
\begin{gather}
  {\bm{u}^{{\rm{ga}}}} \buildrel \Delta \over = {\left[ {M_t^{{\rm{ga,pro}}},M_t^{{\rm{ga,sell}}}} \right]_{t = 1, \ldots ,T}},
\end{gather}
with scenario-dependent decision variables ${[{\bm{u}}_\omega ^{{\rm{ga}}}]_{\omega  \in \Omega }}$.

Because GA is also risk-averse, the CVaR of $-C^{{\rm{ga}}}$ is minimized:
\begin{gather}
  \mathop {\min }\limits_{{{[{\bm{u}}_\omega ^{{\rm{ga}}}]}_{\omega  \in \Omega }}} \;{f^{{\rm{ga}}}}\left( {{{[{\bm{u}}_\omega ^{{\rm{ga}}}]}_{\omega  \in \Omega }},{\bm{v}}} \right) \buildrel \Delta \over = {\rm{CVa}}{{\rm{R}}_\alpha }\left( { - {C^{{\rm{ga}}}}} \right), \label{eq:objcompactga} \\
  \text{s.t.}\   {\bm{h}}_{}^{{\rm{ga}}}\left( {{\bm{u}}_\omega ^{{\rm{ga}}},{\bm{v}}} \right) = \bm{0}, \ \forall \omega \in \Omega :\;{\bm{\lambda} ^{{\rm{ga}}}}, \label{eq:hfunccompactga}  \\
  {\bm{g}}_{}^{{\rm{ga}}}\left( {{\bm{u}}_\omega ^{{\rm{ga}}},{\bm{v}}} \right) \le \bm{0},\ \forall \omega \in \Omega \ : \ {\bm{\mu} ^{{\rm{ga}}}}, \label{eq:gfunccompactga}
\end{gather}
\noindent
where ${f^{{\rm{ga}}}}(\cdot)$ denotes the utility function; ${\bm{h}}^{{\rm{ga}}}(\cdot)$ and ${\bm{g}}^{{\rm{ga}}}(\cdot)$ represent equality and inequality constraints. Because GA production does not depend directly on renewable generation (although it is indirectly affected through market supply-demand and futures trading), the scenario vector $\bm{\phi}_\omega$ does not appear explicitly. After applying the convex linear transformation of CVaR via (\ref{eq:cvarlin})--(\ref{eq:cvarcons}), the problem also reduces to an LP.

\section{Futures-spot dual-market game and equilibrium analysis}
\label{sec:game}

\begin{figure}[t]
	\centering
	\includegraphics[scale=1.05]{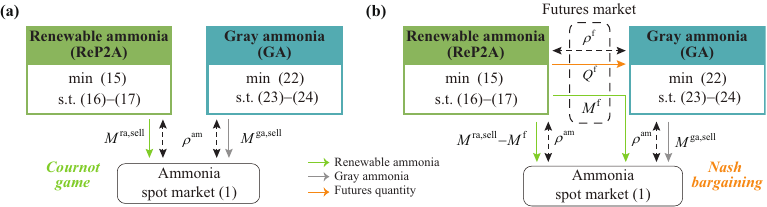}\vspace{-6pt}
	\caption{Structure of the interaction models for ReP2A and GA in (a) the ammonia spot market only (Section \ref{sec:cournot}) and (b) the futures-spot dual market (Section \ref{sec:bargining}). }
	\label{fig:game}
\end{figure}

Based on the operational models of the ReP2A and GA producers in Sections \ref{sec:modelrep2a}--\ref{sec:modelga}, this section develops an interaction model for the ammonia spot and futures markets.
We first analyze the benchmark case without ammonia futures trading. Section \ref{sec:cournot} formulates a Nash-Cournot game for the spot market, as illustrated in Fig.~\ref{fig:game}(a). The resulting equilibrium provides the disagreement point for later bargaining.
Section \ref{sec:bargining} then introduces bilateral ammonia futures trading between ReP2A  and GA producers and develops a Nash bargaining strategy, as shown in Fig.~\ref{fig:game}(b). The mechanism coordinates production and trading decisions to mitigate risks caused by renewable generation uncertainty and improve the utilities of both producers.

\subsection{Nash-Cournot game model of the ammonia spot market}
\label{sec:cournot}

In the benchmark case,  ReP2A and GA sell ammonia simultaneously in the spot market without futures trading. The equilibrium satisfies two conditions:

\emph{1) Under each renewable scenario $\omega \in \Omega$, the ammonia price and total supply satisfy the elasticity relationship (\ref{eq:elasticity});}

\emph{2) The operational decision models of ReP2A and GA, i.e., (\ref{eq:objcompactrep2a})--(\ref{eq:gconsobjcompactrep2a}) and (\ref{eq:objcompactga})--(\ref{eq:gfunccompactga}), both reach optimum.}

One approach to compute the equilibrium is to derive the KKT conditions of the operational optimization problems of ReP2A and GA, convert the complementary slackness conditions into a mixed-integer linear program (MILP), and solve the resulting model \cite{zeng2025planning}. However, the CVaR formulation requires constraints for all scenarios $\omega \in \Omega$, which introduces a large number of binary variables and leads to excessively high computational complexity.

To avoid this issue, we introduce a Gauss-Seidel iterative procedure to compute the equilibrium fixed point. The procedure also reflects the sequential adjustment of market participants.

\emph{Step 1:} Initialize iteration index $k=1$; assign initial spot prices $\rho_{\omega,t}^{{\rm am},(0)}$ for all $\omega \in \Omega$ and $t=1,\ldots,T$; specify the step size $0<\gamma\le1$ and convergence threshold $\epsilon$;

\emph{Step 2:} Treat the current  spot price $\rho_{\omega,t}^{{\rm am},(k-1)}$ as a parameter and solve the ReP2A and GA operational models (\ref{eq:objcompactrep2a})--(\ref{eq:gconsobjcompactrep2a}) and (\ref{eq:objcompactga})--(\ref{eq:gfunccompactga}); obtain optimal decisions and extract ammonia sales $M_{\omega,t}^{{\rm ra,sell},(k)}$ and $M_{\omega,t}^{{\rm ga,sell},(k)}$ for each scenario $\omega$ and period $t$;

\emph{Step 3:} Substitute $M_{\omega,t}^{{\rm ra,sell},(k)}$ and $M_{\omega,t}^{{\rm ga,sell},(k)}$ into the elasticity (\ref{eq:elasticity}) to compute a spot price $\rho_{\omega,t}^{{\rm am}}$; update the price as
$\rho _{\omega ,t}^{{\rm{am,}}\left( k \right)} \leftarrow \left( {1 - \gamma } \right)\rho _{\omega ,t}^{{\rm{am,}}\left( {k - 1} \right)} + \gamma \rho _{\omega ,t}^{{\rm{am}}}$;

\emph{Step 4:} Check convergence. If the changes in $M_{\omega,t}^{{\rm ra,sell}}$, $M_{\omega,t}^{{\rm ga,sell}}$, and $\rho_{\omega,t}^{{\rm am}}$ are below $\epsilon$, terminate and output the equilibrium; otherwise set $k \leftarrow k+1$ and return to \emph{Step 2}.

\begin{remark}
	\label{remark:3}
	The equilibrium obtained here represents the benchmark utilities of the ReP2A and GA without futures trading and serves as the disagreement point for Nash bargaining.
\end{remark}

\subsection{Nash bargaining strategy in the ammonia futures-spot dual market}
\label{sec:bargining}

To ensure that ammonia futures trading benefits both producers, a Nash bargaining framework is adopted. Under Pareto efficiency, the additional utility created by futures trading is allocated between ReP2A and GA through bargaining.

The disagreement point is given by the equilibrium obtained in Section \ref{sec:cournot}. Let the utilities of the ReP2A and GA producers at this point be $f^{{\rm ra*,d}}$ and $f^{{\rm ga*,d}}$. Under a bargaining agreement, the utilities become $f^{{\rm ra*,a}}$ and $f^{{\rm ga*,a}}$. The surplus utilities are therefore ${f^{{\rm{ra*,a}}}} - {f^{{\rm{ra*,d}}}}$ and ${f^{{\rm{ga*,a}}}} - {f^{{\rm{ga*,d}}}}$.
The Nash bargaining solution satisfies four conditions:

\emph{1) Under all renewable scenarios $\omega \in \Omega$, the ammonia price and total supply satisfy the elasticity (\ref{eq:elasticity});}

\emph{2) The operational decision models of ReP2A and GA, i.e., (\ref{eq:objcompactrep2a})--(\ref{eq:gconsobjcompactrep2a}) and (\ref{eq:objcompactga})--(\ref{eq:gfunccompactga}), reach optimum;}

\emph{3) The Nash objective is maximized, i.e., $\text{max}\ F \triangleq \left( {{f^{{\rm{ra*,a}}}} - {f^{{\rm{ra*,d}}}}} \right)\left( {{f^{{\rm{ga*,a}}}} - {f^{{\rm{ga*,d}}}}} \right)$;}

\emph{4) Individual rationality holds, meaning utilities cannot fall below the disagreement, as ${f^{{\rm{ra/ga*,a}}}} \ge {f^{{\rm{ra/ga*,d}}}}$.}

To obtain the bargaining outcome, we construct an iterative strategy in which both producers update futures trading and operational decisions. To maximize $F$, the gradients of the Nash objective with respect to futures trading position $\bm{Q}^{\rm f}$ and price $\bm{\rho}^{\rm f}$ are	derived:
\begin{gather}
	{{\Delta F} \over {\Delta \bm{Q}_{}^{\rm{f}}}} = \left( {{f^{{\rm{ra}}}} - {f^{{\rm{ra*,d}}}}} \right){{\Delta {f^{{\rm{ra*}}}}} \over {\Delta {\bm{Q}}_{}^{\rm{f}}}} + \left( {{f^{{\rm{ga}}}} - {f^{{\rm{ga*,d}}}}} \right){{\Delta {f^{{\rm{ga*}}}}} \over {\Delta {\bm{Q}}_{}^{\rm{f}}}}, \label{eq:gradquan} \\
	{{\Delta F} \over {\Delta \bm{\rho}_{}^{\rm{f}}}} = \left( {{f^{{\rm{ra}}}} - {f^{{\rm{ra*,d}}}}} \right){{\Delta {f^{{\rm{ra*}}}}} \over {\Delta \bm{\rho} _{}^{\rm{f}}}} + \left( {{f^{{\rm{ga}}}} - {f^{{\rm{ga*,d}}}}} \right){{\Delta {f^{{\rm{ga*}}}}} \over {\Delta \bm{\rho} _{}^{\rm{f}}}} \label{eq:gradpr},
\end{gather}
\noindent
where sensitivities of the optimal utilities ${f^{{\rm{ra/ga*}}}}$ w.r.t. ${\bm{Q}}^{\rm f}$ and $\bm{\rho}^{\rm f}$ are obtained from the Lagrange multipliers of the operational models (\ref{eq:objcompactrep2a})--(\ref{eq:gconsobjcompactrep2a}) and (\ref{eq:objcompactga})--(\ref{eq:gfunccompactga}), as
\begin{gather}
	{{\Delta {f^{{\rm{ra/ga*}}}}} \over {\Delta {\bm{Q}}_{}^{\rm{f}}}} = \left[ {{{\partial {{\bm{h}}^{{\rm{ra/ga}}}}} \over {\partial {\bm{Q}}_{}^{\rm{f}}}}} \right]{\left( {{\bm{\lambda} ^{{\rm{ra/ga}}}}} \right)^{\rm{T}}} + \left[ {{{\partial {{\bm{g}}^{{\rm{ra/ga}}}}} \over {\partial {\bm{Q}}_{}^{\rm{f}}}}} \right]{\left( {{\bm{\mu} ^{{\rm{ra/ga}}}}} \right)^{\rm{T}}}. \label{eq:sen}
\end{gather}

The bargaining procedure combines the Gauss-Seidel updates used in Section \ref{sec:cournot} with gradient descent, as follows:

\emph{Step 1:} Solve the Nash-Cournot equilibrium in Section \ref{sec:cournot} to obtain the disagreement utilities $f^{{\rm ra*,d}}$ and $f^{{\rm ga*,d}}$; use the resulting spot price as the initial value $\rho_{\omega,t}^{{\rm am},(0)}$; specify step sizes $\gamma$, $\beta^\rho$, $\beta^Q$ and convergence threshold $\epsilon$;

\emph{Step 2:} Initialize iteration index $k=1$ and futures trading variables $\bm{Q}^{\rm f,(0)}$ and $\bm{\rho}^{\rm f,(0)}$ for each scenario $\omega \in \Omega$ and each period $t=1,\ldots, T$;

\emph{Step 3:} Given $\rho_{\omega,t}^{{\rm am},(k-1)}$, $\bm{Q}^{\rm f,(k-1)}$, and $\bm{\rho}^{\rm f,(k-1)}$, solve the ReP2A and GA operational models (\ref{eq:objcompactrep2a})--(\ref{eq:gconsobjcompactrep2a}) and (\ref{eq:objcompactga})--(\ref{eq:gfunccompactga}) to obtain $M_{\omega,t}^{{\rm ra,sell},(k)}$ and $M_{\omega,t}^{{\rm ga,sell},(k)}$ for all periods and scenarios;

\emph{Step 4:} Update the spot price using the elasticity (\ref{eq:elasticity}) and the step-size rule
$\rho _{\omega ,t}^{{\rm{am,}}\left( k \right)} \leftarrow \left( {1 - \gamma } \right)\rho _{\omega ,t}^{{\rm{am,}}\left( {k - 1} \right)} + \gamma \rho _{\omega ,t}^{{\rm{am}}}$;

\emph{Step 5:} Update the futures price and quantity using gradient descent, as
$\bm{\rho} _{}^{{\rm{f,}}\left( k \right)} \leftarrow \bm{\rho} _{}^{{\rm{f,}}\left( {k - 1} \right)} - {\beta ^\rho } \dfrac{\Delta F} {\Delta \bm{\rho} _{}^{\rm{f}}}$,
$\bm{Q} _{}^{{\rm{f,}}\left( k \right)} \leftarrow \bm{Q} _{}^{{\rm{f,}}\left( {k - 1} \right)} - {\beta ^Q } \dfrac{\Delta F} {\Delta \bm{Q} _{}^{\rm{f}}}$;

\emph{Step 6:} Check convergence. If the changes in $M_{\omega,t}^{{\rm ra,sell}}$, $M_{\omega,t}^{{\rm ga,sell}}$, $\rho_{\omega,t}^{{\rm am}}$, ${\bm{Q}}^{\rm f}$, and $\bm{\rho}^{\rm f}$ are below $\epsilon$, terminate; otherwise set $k \leftarrow k+1$ and return to \emph{Step 3};

\emph{Step 7:} Verify the solution. If the utilities of both producers exceed the disagreement point, the result is accepted as the Nash bargaining agreement; otherwise the bargaining fails.

\begin{remark}
	\label{remark:4}
	An alternative formulation is a non-cooperative generalized Nash equilibrium problem (GNEP), as in \cite{zeng2025planning}. However, numerical experiments indicate that such equilibria may not exist or may produce utilities below the disagreement level. This study therefore adopts the Nash bargaining framework, which ensures mutually beneficial outcomes.
\end{remark}

\section{Case studies}
\label{sec:case}

\subsection{Case settings}
\label{sec:setting}

The case study is based on a ReP2A demonstration project in northern China \cite{zhu2026exploring,zeng2025optimal}. The system includes a 450 MW wind farm that supplies electricity to an ammonia plant with an annual capacity of 200,000 t (22.83 t/h). Each MWh of renewable electricity produces ${\eta ^{{\rm{p2a}}}} = 0.1030$ t/MWh of ammonia. The fixed operating cost of the ReP2A system, including operation and maintenance of the wind farm, transmission facilities, hydrogen production, and the ammonia plant, is 42,000 CNY/h.
A conventional GA plant is assumed to have the same annual capacity of 200,000 t. Its fixed operating cost is 41,200 CNY/h and its variable production cost is 1,320 CNY/t.

The ammonia demand is defined such that the market price reaches ${\rho ^{{\rm{max}}}} = 4,850$ CNY/t when supply is zero. The elasticity coefficient is ${k^{{\rm{am}}}} = 16.44$ t$^2$/CNY. Under this specification, when monthly supply reaches 100,000 t the price decreases by 506.9 CNY/t. The parameters are summarized in Table \ref{tab:para}.

\begin{table}[b]\scriptsize\centering
	\renewcommand{\arraystretch}{1.15}
	\caption{Base-case parameters of the case study}\vspace{6pt}
	\label{tab:para}
	\begin{tabular}{ ccc }\hline\hline
		Parameter &  Meaning  &  Value \\ \hline
		${\overline{M}^{{\rm{ra,pro}}}}$,  ${\underline{M}^{{\rm{ra,pro}}}}$     & ReP2A hourly production upper/lower bounds & 22.83 t/h, 0    \\
		${\overline{M}^{{\rm{ga,pro}}}}$,  ${\underline{M}^{{\rm{ga,pro}}}}$       & GA hourly production upper/lower bounds  & 22.83 t/h, 0 	  \\
		${\eta ^{{\rm{p2a}}}}$			& ReP2A energy conversion coefficient	& 0.1030 t/MWh 	\\
		${c^{{\rm{ra}}}}$  				&  ReP2A hourly production cost			& 42,000 CNY   	\\
		$c_{\rm{0}}^{{\rm{ga}}}$ 		& GA hourly fixed cost					& 41,000 CNY	\\
		$c_1^{{\rm{ga}}}$				& GA variable cost per ton of ammonia	& 1,320 CNY/t	\\
		${\rho ^{{\rm{max}}}}$, ${k^{{\rm{am}}}}$	& Ammonia market price coefficients & 4,850 CNY/t, 16.44 t$^2$/CNY \\
		$\alpha$						& CVaR confidence level					& 0.5\\ \hline\hline
	\end{tabular}
\end{table}

\begin{figure}[t]
	\centering
	\includegraphics[scale=0.9]{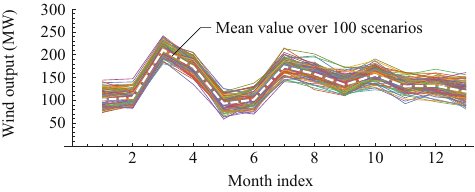}\vspace{-6pt}
	\caption{Scenario set of monthly average wind power output used in the base case.}
	\label{fig:scen}
\end{figure}

\begin{figure}[t]
	\centering
	\includegraphics[scale=0.9]{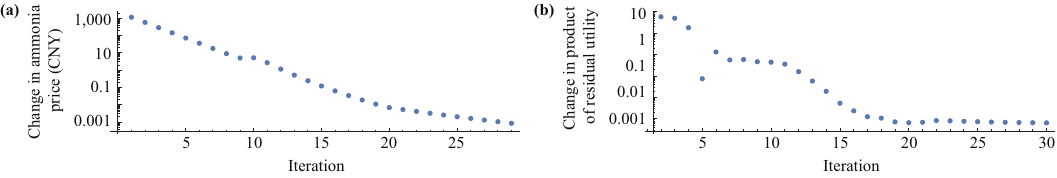}\vspace{-6pt}
	\caption{(a) Convergence of ammonia price solving Nash-Cournot equilibrium under the spot-only market. (b) Convergence process of the futures-spot joint market equilibrium.}
	\label{fig:conv2}
\end{figure}

Monthly average wind generation follows measured data from the project site \cite{zhu2026exploring,zeng2025optimal}, as shown in Fig. \ref{fig:scen}. Uncertainty is modeled by adding a uniformly distributed disturbance within $\pm10\%$ of the rated wind capacity, i.e., $[-45,45]$ MW, followed by a moving-average filter with window width 2. A total of $N^{\text{s}} = 100$ scenarios are generated to represent renewable variability.
Risk preferences are represented using CVaR with confidence level $\alpha = 0.5$, corresponding to moderate risk aversion. The influence of the confidence level is analyzed later in Section \ref{sec:sen}. All models are implemented in \emph{Wolfram Mathematica 13.0}, and the production optimization problems are solved using \emph{Mosek 11.1}.

\subsection{Base case analysis}
\label{sec:base}

\subsubsection{Spot market equilibrium without futures trading}
\label{sec:casecournot}

We first consider the benchmark case without futures trading, in which the ReP2A and GA producers compete only in the ammonia spot market. The resulting equilibrium serves as the disagreement point in the Nash bargaining analysis and provides a baseline for evaluating the proposed futures mechanism.

\begin{figure}[t]
	\centering
	\includegraphics[scale=0.9]{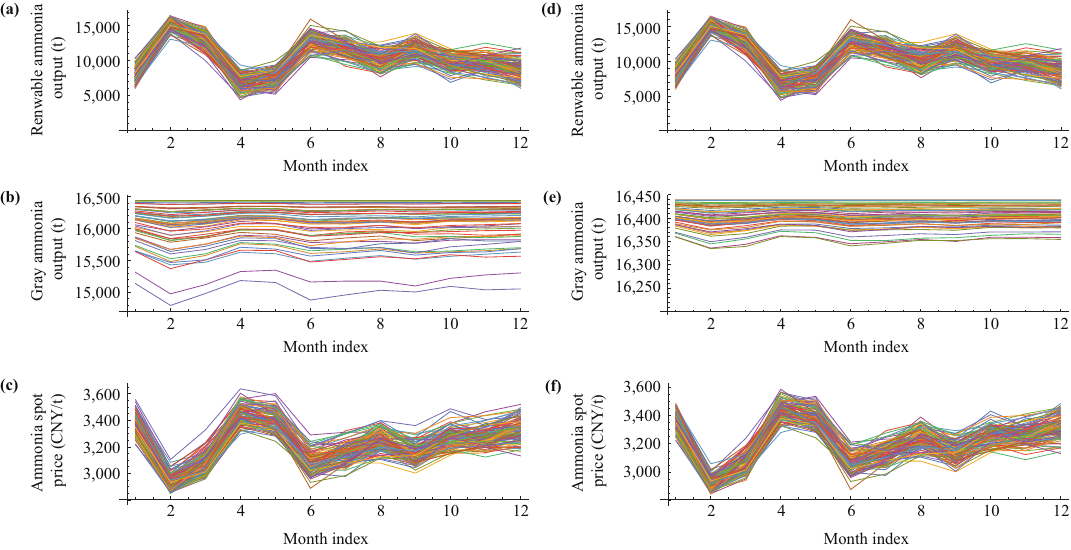}\vspace{-6pt}
	\caption{Operational results of under spot market only (a) ReP2A production; (b) GA production; (c) ammonia spot market price; and after introducing the futures mechanism (d) ReP2A production; (e) GA production; (f) ammonia spot market price, across 100 scenarios.
	}
	\label{fig:op2}
\end{figure}

Using the algorithm described in Section \ref{sec:cournot}, the Nash-Cournot equilibrium is obtained. The convergence of the ammonia price during iterations is shown in Fig. \ref{fig:conv2}(a). The Gauss-Seidel algorithm reduces the price change to below $0.001$ CNY/t within approximately 20 iterations.
At equilibrium, the utilities of ReP2A and GA are $2.963 \times 10^{7}$ CNY and $1.495 \times 10^{7}$ CNY. These values serve as the disagreement utilities in the bargaining analysis in Section \ref{sec:casebargining}.

The resulting production and spot prices are shown in Fig. \ref{fig:op2}(a)-- \ref{fig:op2}(c). Due to renewable uncertainty, ReP2A output varies between 5,000 t and 16,000 t across months and scenarios. This variability propagates to total market supply and causes the ammonia price to fluctuate between 2,800 CNY/t and 3,600 CNY/t. Such price volatility affects the GA's decisions. When the market price approaches production cost, the GA producer reduces output to avoid losses, lowering production by up to 1,500 t from its technical maximum. Consequently, GA production ranges between 14,800 t and 16,400 t.

\subsubsection{Effect of introducing renewable ammonia futures}
\label{sec:casebargining}

The renewable ammonia futures mechanism described in Section \ref{sec:futures} is then introduced. The equilibrium of the coupled spot-futures market is solved using the hybrid  algorithm proposed in Section \ref{sec:bargining}. For \emph{Mode 2 (conventional commodity-futures type)}, Nash bargaining fails because no futures transaction allows both producers to achieve utilities above the disagreement point. This confirms the theoretical analysis in Section \ref{sec:futures} and indicates that this mode cannot effectively hedge renewable-induced risks. The following analysis therefore focuses on \emph{Mode 1 (wind-futures type)}.

For Mode 1, the convergence of the Nash objective (product of surplus utilities of ReP2A and GA) is shown in Fig. \ref{fig:conv2}(b). The algorithm converges within about 20 iterations.
Table \ref{tab:comp} compares the outcomes with and without futures trading. After introducing renewable ammonia futures, the utility of ReP2A increases to $3.114 \times 10^{7}$ CNY, while the GA's utility rises to $1.647 \times 10^{7}$ CNY, with improvements of 5.075\% and 10.14\%, respectively. The combined utility increases by 7.063\%, indicating improved overall welfare.


\begin{figure}[t]
	\centering
	\includegraphics[scale=0.9]{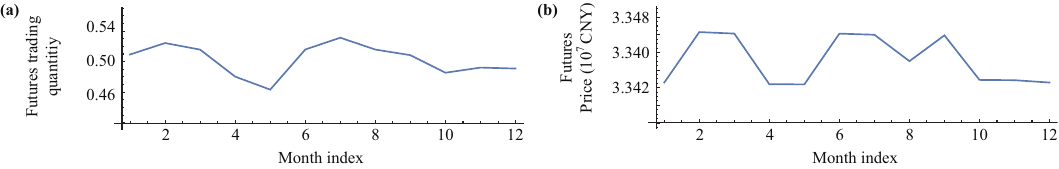}\vspace{-6pt}
	\caption{Monthly renewable ammonia futures trading in the base case. (a) Futures position. (b) Futures price.}
	\label{fig:fu}
\end{figure}

\begin{table}[tb]\footnotesize\centering
	\renewcommand{\arraystretch}{1.4}
	\caption{Comparison of utilities before and after introducing renewable ammonia futures.}\vspace{6pt}
	\label{tab:comp}
	\begin{tabular}{ cccc }\hline\hline
		Market model 		& ReP2A utility (CNY)  & GA utility (CNY) &  Total utilities (CNY) \\ \hline
		w/o futures trading & $2.963 \times 10^7$	& $1.495 \times 10^7$	& $4.446 \times 10^7$	\\		
		w/ futures trading & \tabincell{c}{ $3.113 \times 10^7$ \\ {(\bf+5.075\%})}	& \tabincell{c}{ $1.647 \times 10^7$ \\ {(\bf+10.14\%})}	& \tabincell{c}{ $4.760 \times 10^7$ \\ {(\bf+7.063\%})}	\\ \hline\hline
	\end{tabular}
\end{table}

\begin{figure}[tb]
	\centering
	\includegraphics[scale=0.95]{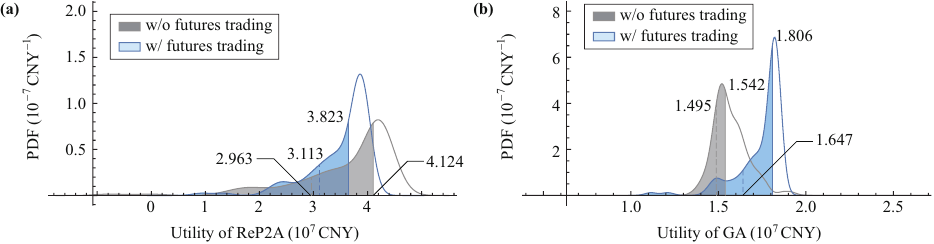}\vspace{-6pt}
	\caption{Probability density functions of utilities for ReP2A and GA (a) before, and (b) after introducing the futures mechanism.}
	\label{fig:pdfra}
\end{figure}

Production outcomes and spot prices with futures trading are shown in Fig. \ref{fig:op2}(d)--\ref{fig:op2}(f), while futures trading positions and prices are presented in Fig. \ref{fig:fu}. ReP2A output remains unchanged because it is constrained by renewable generation. In contrast, GA production increases and becomes more stable, with fluctuations shrinking from $[14,800,16,400]$ t to $[16,300,16,400]$ t. The increase in GA output slightly raises total supply and reduces market prices, benefiting consumers.

Fig. \ref{fig:fu} shows that both the futures position and price vary only slightly during the year. The trading quantity remains between 0.47 and 0.53, while the price stays around $3.340 \times 10^{7}$ CNY with fluctuations within $\pm0.15\%$. These patterns follow the seasonal variation of renewable generation shown in Fig. \ref{fig:scen}. The trading positions close to 50\% reflect the symmetric production capacities of the two producers in the base case. Under the Nash bargaining solution, the production risk associated with renewable uncertainty is shared approximately equally between them.

To further illustrate the risk-hedging effect, Fig. \ref{fig:pdfra} presents the probability density functions (PDFs) of utilities before and after introducing futures trading. For ReP2A, futures contracts convert part of the uncertain spot revenue into stable futures income, narrowing the utility distribution and reducing downside risk. The probability of outcomes below $2 \times 10^{7}$ CNY declines significantly. Although the VaR at $\alpha = 0.5$ decreases from $4.124 \times 10^{7}$ CNY to $3.823 \times 10^{7}$ CNY, the CVaR increases from $2.963 \times 10^{7}$ CNY to $3.113 \times 10^{7}$ CNY.
For GA, futures trading allows part of the RA output to be purchased at a stable price. This shifts the utility distribution to the right. As a result, the VaR at $\alpha = 0.5$ increases from $1.542 \times 10^{7}$ CNY to $1.806 \times 10^{7}$ CNY, and the CVaR increases from $1.495 \times 10^{7}$ CNY to $1.647 \times 10^{7}$ CNY.

\subsection{Sensitivity analysis and discussion}
\label{sec:sen}

Additional sensitivity analyses examine how renewable uncertainty, CVaR confidence levels, and production capacities influence the effectiveness of the proposed futures mechanism.

\subsubsection{Renewable generation uncertainty}
\label{sec:senuncertainty}

\begin{figure}[t]
	\centering
	\includegraphics[scale=0.9]{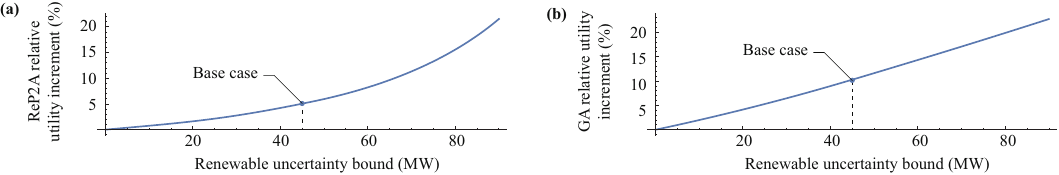}\vspace{-6pt}
	\caption{Relative utility improvement from the futures mechanism under different levels of renewable generation uncertainty. (a) ReP2A. (b) GA. }
	\label{fig:senuncertain}
	\vspace{12pt}
	\centering
	\includegraphics[scale=0.9]{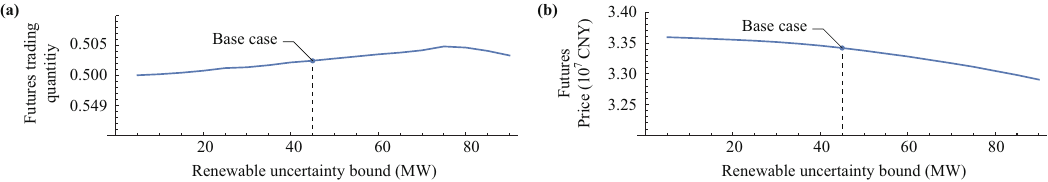}\vspace{-6pt}
	\caption{Average futures trading under different levels of renewable generation uncertainty. (a) Futures position. (b) Futures price.}
	\label{fig:senuncertainfu}
\end{figure}	

The disturbance applied to renewable output is expanded from $[0,0]$ to $[-90,90]$ MW to represent increasing uncertainty. Other parameters remain the same as in Section \ref{sec:setting}.
Fig. \ref{fig:senuncertain} shows the resulting utility improvements from futures trading. When no uncertainty exists, futures trading provides no benefit. As uncertainty increases, utility gains become more significant for both sides. The improvement grows faster for ReP2A, indicating that futures trading mitigates production risk more than market price risk.

Fig. \ref{fig:senuncertainfu} shows the corresponding futures trading outcomes. The trading position remains close to 0.5 and the price stays near $3.350 \times 10^{7}$ CNY. As uncertainty increases, the futures position rises slightly while the price declines marginally, reflecting the higher risk of ReP2A production. Overall, however, the trading outcomes remain close to the base case in Section \ref{sec:casebargining}.

\subsubsection{CVaR confidence level}
\label{sec:senalpha}

\begin{figure}[t]
	\centering
	\includegraphics[scale=0.9]{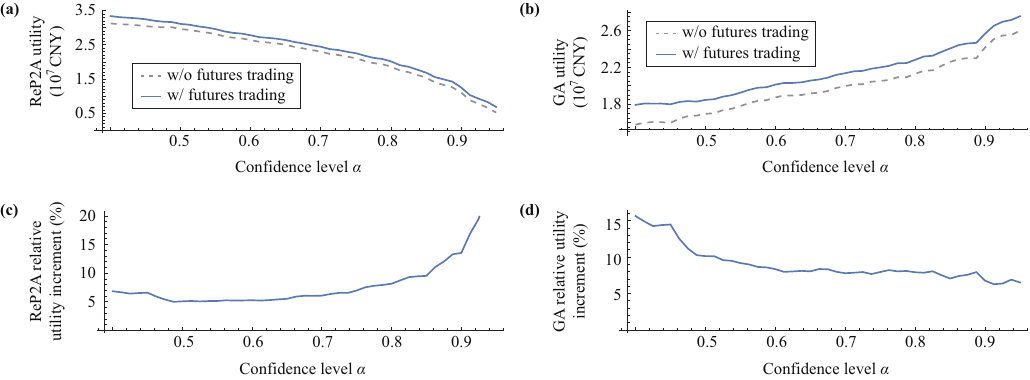}\vspace{-6pt}
	\caption{Impact of CVaR confidence levels on the utilities of ReP2A and GA. (a) ReP2A utility. (b) GA utility. (c) Relative utility improvement for ReP2A from the futures mechanism. (d) Relative utility improvement for GA from the futures mechanism. }
	\label{fig:alpha}
	\centering
	\vspace{12pt}
	\includegraphics[scale=0.9]{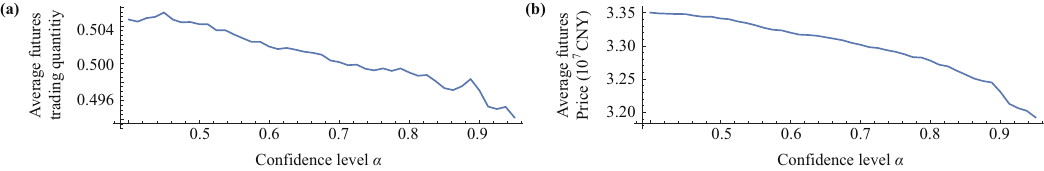}\vspace{-6pt}
	\caption{Futures trading under different CVaR confidence levels. (a) Futures position. (b) Futures price. }
	\label{fig:alphafu}
\end{figure}

The CVaR confidence level $\alpha$ determines how strongly extreme downside outcomes influence decision making. To analyze its influence, $\alpha$ is varied from 0.4 to 0.95.
Fig. \ref{fig:alpha} shows that higher confidence levels increase the benefit of futures trading for ReP2A. The improvement becomes particularly strong when $\alpha > 0.8$ because futures trading substantially reduces the left tail of the utility distribution as shown in Fig. \ref{fig:pdfra}.
For the GA producer, the benefit mainly arises from the rightward shift of its utility distribution. However, this shift slightly widens the distribution, leading to minor deterioration in extreme outcomes. Consequently, the relative improvement decreases slightly as $\alpha$ increases.

Fig. \ref{fig:alphafu} shows that the futures position remains close to 0.5 across all confidence levels, indicating balanced risk sharing. The futures price decreases slightly (approx. 5\%) as $\alpha$ increases because stronger aversion to extreme losses reduces the willingness to bear risk.

When $\alpha = 0$, both producers become risk-neutral and maximize expected utility. In this case, the trading equilibrium depends on the initial futures positions in the algorithm, but the utilities remain identical to those without futures trading. This confirms that the proposed mechanism primarily functions as a risk-hedging instrument.

\subsubsection{Capacity scales of ReP2A and GA}
\label{sec:sensize}

\begin{figure}[tb]
	\centering
	\includegraphics[scale=0.95]{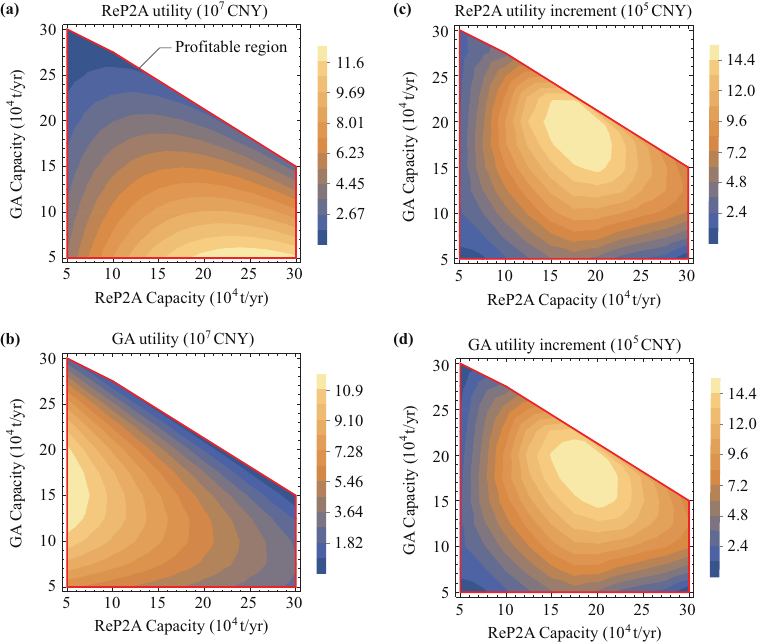}\vspace{-6pt}
	\caption{The effects of ReP2A-GA capacity combination. (a) ReP2A utility without futures trading. (b) GA utility without futures trading. (c) ReP2A utility improvement after introducing futures trading. (d) GA utility improvement after introducing futures trading. }
	\label{fig:sensize}
\end{figure}

The impact of production capacity is examined by varying the capacities of both producers from 50,000 to 300,000 t/year. Fixed costs are scaled proportionally with capacity while other parameters remain unchanged. Figs. \ref{fig:sensize}(a) and \ref{fig:sensize}(b) show the utilities without futures trading. The feasible region where both utilities are positive is highlighted by a red box.

For both producers, utility initially increases with capacity because higher output raises sales. Beyond a certain level, however, market saturation reduces prices through the elasticity relation (\ref{eq:elasticity}), which lowers profits. Increasing the competitor's capacity further intensifies competition and reduces utilities.

Figs. \ref{fig:sensize}(c) and \ref{fig:sensize}(d) show the utility improvements after introducing futures trading. Benefits appear throughout the feasible region, with the largest gains occurring when both capacities lie between 150,000 and 200,000 t/year. In this range, market competition and production risk are both significant, making risk sharing particularly valuable. When the capacity difference becomes large, the effectiveness of futures trading declines because the risk exposures of the two producers become less balanced.

\begin{figure}[t]
	\centering
	\includegraphics[scale=0.95]{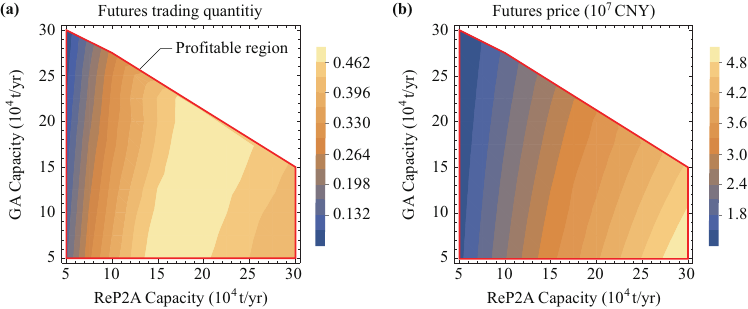}\vspace{-6pt}
	\caption{Futures trading under different ReP2A-GA capacity combinations. (a) Futures position. (b) Futures price.}
	\label{fig:sizefu}
\end{figure}

Fig. \ref{fig:sizefu} shows the corresponding futures trading outcomes. Trading volume is highest when ReP2A capacity lies between 150,000 and 200,000 t/year. When capacity is smaller, price risk declines and trading volume falls. When capacity becomes very large, market saturation reduces the value of transferring production through futures contracts.
Another interesting thing is that reducing GA capacity does not remarkably reduce trading volume. Even when GA capacity is only 50,000 t/year, futures trading remains substantial. This suggests that institutions could potentially participate in the market through futures trading alone. This possibility is examined in detail in Section \ref{sec:nptp}.

\subsubsection{Discussion on non-producing trading participants (NPTP)}
\label{sec:nptp}

The previous results indicate that a participant without ammonia production could potentially engage in renewable ammonia futures trading. To examine this possibility, a non-producing trading participant (NPTP) is introduced.
The GA model from Section \ref{sec:modelga} is used, but production capacity and costs are set to zero, leaving only the trading function. Because financial institutions typically exhibit lower risk aversion, the CVaR confidence level for the NPTP is set to $\alpha = 0.2$, while the ReP2A producer retains $\alpha = 0.5$.

\begin{figure}[tb]
	\centering
	\includegraphics[scale=0.95]{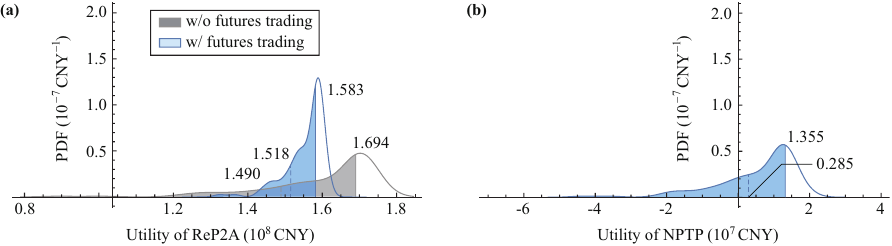}\vspace{-6pt}
	\caption{Probability density functions of utilities for ReP2A and NPTP (a) before and (b) after introducing the futures mechanism. }
	\label{fig:nonprod}
\end{figure}

Fig. \ref{fig:nonprod} shows the resulting utility distributions. Without futures trading, the ReP2A utility is $1.490 \times 10^{8}$ CNY. After introducing futures trading, it increases to $1.518 \times 10^{8}$ CNY, corresponding to a 1.915\% improvement. The NPTP obtains zero utility without trading but earns $2.846 \times 10^{6}$ CNY by participating in the futures market.
This outcome occurs because the less risk-averse NPTP is willing to assume part of the renewable production risk in exchange for higher expected returns. The average futures position increases to 0.69, indicating that the NPTP absorbs a larger share of the risk.
However, if the NPTP becomes more risk-averse than the ReP2A producer (i.e., its $\alpha$ exceeds that of ReP2A), Nash bargaining fails and no futures trading occurs. This confirms that the mechanism relies on differences in risk preferences to enable effective risk sharing.

These results suggest an important policy implication. In markets dominated by ReP2A producers, regulators could allow financial institutions to participate as NPTPs. Such participants could help hedge renewable production risk and support investment in green ammonia industries.

%

\section{Conclusions}
\label{sec:conclusion}

Variability in renewables leads to fluctuations in ammonia output in ReP2A systems, which propagate to market supply and prices, creating revenue risks for both ReP2A and conventional GA producers. To address this challenge, this paper proposes a \emph{renewable ammonia futures} mechanism. A coupled spot-futures market game model with a Nash bargaining framework is developed to coordinate production and trading decisions. Case studies based on a real-life project demonstrate the effectiveness of the proposed approach. The main findings include:

\begin{enumerate}
	\item Renewable generation uncertainty affects ammonia prices through fluctuations in production, reducing revenue stability for both ReP2A and GA producers. Renewable ammonia futures allow part of the uncertain production capacity to be converted into predetermined revenue, thereby reducing downside risk. In the base case, the CVaR utilities of ReP2A and GA increase by approximately 5.103\% and 10.14\%, respectively.

	\item The effectiveness of the futures mechanism depends on renewable uncertainty and risk preferences. As renewable output variability increases, the benefits of futures trading become more pronounced for both producers, particularly for ReP2A. Higher CVaR confidence levels further reduce downside risk for ReP2A, while the relative benefit for GA decreases slightly.

	\item The impact of the mechanism also depends on the capacity configuration of the two producers. The largest utility improvements occur when ReP2A and GA have comparable capacities and market competition is intense. When the capacity difference becomes large, the effectiveness of futures trading declines because risk sharing becomes less balanced.
	
	\item Even in markets without GA producers, non-producing trading participants can profit by assuming part of the renewable production risk while improving the revenue stability of ReP2A producers. This finding highlights the potential role of financial institutions in supporting risk management in renewable ammonia markets.
\end{enumerate}

Future research could further consider technological developments in ReP2A systems, hydrogen and ammonia storage options, evolving market structures, demand uncertainty, and long-term investment to better assess the role of integrated financial-engineering mechanisms in supporting the green transition of the energy and chemical industries.

\section*{Acknowledgement}

The authors gratefully acknowledge the financial support from the National Key Research and Development Program of China (2021YFB4000503) and the National Natural Science Foundation of China (52377116).

\section*{Declaration of Interest}

None.

\section*{Data Availability}

The data related to this work are available upon request.

%


\end{document}